\documentclass[final]{siamltex}

\usepackage{graphicx}
\usepackage{amssymb,amsmath}
\usepackage{comment}
\usepackage{url} 
\usepackage{hyperref}
\usepackage{subcaption}
\usepackage{algorithm}
\usepackage{algpseudocode}

\newcommand{\Input}{\textbf{Input: }}

\newtheorem{estimate}[theorem]{Estimate}

\usepackage{bbm,color}

\def\R{\mathbbm{R}}
\def\deg{d}         
\def\mrpoly{\pi}    

\def\pof{\mbox{{\it pof}}@}

\mathcode`\@="8000
{\catcode`\@=\active \gdef@{\mkern1mu}} 

\begin{document}

\title{Toward Efficient Polynomial Preconditioning for GMRES
\footnotemark[1] }

\author{Jennifer A. Loe\footnotemark[2]
\and Ronald B. Morgan\footnotemark[3] }

\maketitle

\renewcommand{\thefootnote}{\fnsymbol{footnote}}
\footnotetext[1]{The second author was supported by NSF grant DMS-1418677.  He is the corresponding author.}
\footnotetext[2]{Center for Computing Research, Sandia National Laboratories, Albuquerque, NM, 87123 ({\tt jloe@sandia.gov}).  Sandia National Laboratories is a multimission laboratory managed and operated by National Technology and Engineering Solutions of Sandia, LLC, a wholly owned subsidiary of Honeywell International, Inc., for the U.S. Department of Energy's National Nuclear Security Administration under contract DE-NA-0003525. This paper describes objective technical results and analysis. Any subjective views or opinions that might be expressed in the paper do not necessarily represent the views of the U.S. Department of Energy or the United States Government. SAND2022-0307 O}
\footnotetext[3]{Department of Mathematics, Baylor
University, Waco, TX 76798-7328 ({\tt Ronald\_Morgan@baylor.edu}).}
\renewcommand{\thefootnote}{\arabic{footnote}}

\begin{abstract}
We present a polynomial preconditioner for solving large systems of linear equations.  The polynomial is derived from the minimum residual polynomial (the GMRES polynomial) and is more straightforward to compute and implement than many previous polynomial preconditioners.  Our current implementation of this polynomial using its roots is naturally more stable than previous methods of computing the same polynomial. We implement further stability control using added roots, and this allows for high degree polynomials.  We discuss the effectiveness and challenges of root-adding and give an additional check for stability.  In this paper, we study the polynomial preconditioner applied to GMRES; however it could be used with any Krylov solver.  This polynomial preconditioning algorithm can dramatically improve convergence for some problems, especially for difficult problems, and can reduce dot products by an even greater margin.  
\end{abstract}



\pagestyle{myheadings}
\thispagestyle{plain}
\markboth{J. A. LOE and R. B. MORGAN}{NEW POLYNOMIAL PRECONDITIONED GMRES}

\section{Introduction}
\label{Sec:Intro}
We present a new approach for polynomial preconditioning the Krylov method GMRES~\cite{SaSc} for solving $Ax = b$.  While we assume that the matrix $A$ is large and nonsymmetric, the same approach could also be applied effectively in the symmetric case.  We also assume that $A$ is real-valued; minor modifications allow using the same polynomials for a complex-valued matrix.  As in~\cite{PPG},  we use the minimum residual polynomial.  This is derived from a preliminary GMRES run and is also called the GMRES polynomial. So GMRES is used in two ways: initially it finds the polynomial; then restarted GMRES with the polynomial preconditioner solves the linear equations.  The GMRES polynomial forms a general-purpose preconditioner that is compatible with any system of linear equations and can be composed with any standard preconditioner.  Our new implementation gives a significant advancement over the method in~\cite{PPG}, because it is naturally more stable for high degree polynomials.  And with a high degree polynomial, there is more potential for it to accurately approximate $A^{-1}$.  Additional stability modifications via added roots allow us to stably use polynomials of even higher degrees.  The algorithm to generate the polynomial is simple to implement, and stability modifications can be automated.  This new polynomial preconditioner can significantly improve convergence of GMRES for difficult problems and can also greatly reduce orthogonalization expenses.  

We refer to polynomial preconditioning as being unstable when it does not give accurate results for the final residual norm $\|b-A x\|$, where $x$ is the solution generated by the polynomial preconditioned GMRES.  There can be different reasons for the instability.  For the power basis method in~\cite{PPG}, coefficients for the terms of the polynomial $p$ (given
in the next paragraph) may be inaccurately computed because the power basis becomes ill-conditioned.  The polynomial $p$ is computed differently in~\cite{NaReTr} (coefficients of $p$ are found from the Arnoldi iteration), but as pointed out in~\cite{PPG}, it can also be unstable.  In this paper, the inaccuracy is not because of a poorly computed polynomial, but instead due to the ill-conditioned nature of the polynomial itself.  High degree GMRES polynomials can have very steep slopes and when this happens at an eigenvalue, it can cause inaccuracy in computations.  We still refer to this situation as being an unstable method rather than being badly conditioned, because the use of such a polynomial causes the method to be inaccurate.  Our solution to this problem is not to attempt more accurate computation with this ill-conditioned polynomial, but instead to adjust the polynomial to be better conditioned.  In this paper, we do not give theoretical results about this stability, but do illustrate with computations. 

With a polynomial $p$ and right preconditioning, the linear equations problem becomes 
\begin{align}
Ap(A)y &= b, \label{eqn:ppsys} \\  
x &= p(A)y.  \label{eqn:ppsys2}
\end{align}
The rewritten system of linear equations uses the matrix $Ap(A)$, which typically has an improved spectrum for solving with GMRES.  
We let $\phi(\alpha) \equiv \alpha p(\alpha)$ be the polynomial corresponding to $Ap(A)$, where $\alpha$ is an independent variable for the polynomial.  Then equation (\ref{eqn:ppsys}) becomes 
\[\phi(A)y = b.\]  
We let the degree of $p$ be $d-1$ so $\phi$ has degree $d$.  

All Krylov methods find approximate solutions that involve polynomials of the matrix $A$.  However GMRES is limited to polynomials of fairly low degree if it is restarted.  Polynomial preconditioning allows the solution obtained from restarted GMRES to be represented by high-degree polynomials.  Suppose that GMRES restarts after $m$ iterations.  Then GMRES builds the Krylov subspace $Span\{r, Ar, A^2r,$ $\ldots $  $A^{m-1}r\}$, where $r$ is both the residual vector from the previous cycle and the right-hand side of the current system of equations.  An approximate solution from this subspace can be written as $\omega(A)r$, where $\omega$ is a polynomial of degree $m-1$.  In contrast, polynomial preconditioned GMRES has the subspace $Span\{r, \phi(A)r, (\phi(A))^2 r, \ldots ,$ $(\phi(A))^{m-1} r\}$, so an approximate solution for $\phi(A)y = r$ can be written as $y = \omega(\phi(A))r$.  The composite polynomial $\omega \circ \phi$ has degree $(m-1)*d$.  GMRES now can have high-degree polynomials even though it only builds a subspace of dimension $m$.  These high degree polynomials can give much faster convergence for some problems, as will be demonstrated with examples.  

Another way to explain the effectiveness of polynomial preconditioning is that the spectrum with $\phi(A)$ is generally improved compared to the original spectrum of $A$.  The GMRES polynomial $\pi(\alpha)$ is one at the origin and tries to be near zero at the eigenvalues of $A$, so $\phi(\alpha) = 1 - \pi(\alpha)$ is zero at the origin and near to one over the spectrum of $A$.  (See Figure \ref{fig:MembplusGoodvsBadPoly} (a) for a $\phi$ polynomial graph.)  If effective, the polynomial maps most of the eigenvalues near to one and spaces out the smallest eigenvalues, giving an easier problem for a Krylov method.  

Many polynomial preconditioners have been proposed; see for example~\cite{La52B,Sti58,Rutis,Sa84b,Sa87b,As87,SmSa,FiRe,AsMaOt,Jo,vGi95,Sa03,Th06,seed,LiXiVeYaSa,PPG}. However, they often are complicated to implement for nonsymmetric problems and may require computing estimates to bound a complex spectrum. 
These polynomials can also be unstable and thus have significant computational error even when the linear equations are reasonably well-conditioned~\cite{LiWeSzStabPP,PPG}.  Thus, polynomial preconditioners are not typically used in practice.  Recent work indicates renewed interest in polynomial preconditioning, but this work has been mostly for the symmetric case (e.g.~\cite{LiXiVeYaSa}).  This renewed interest likely arises from the need for scalable preconditioners for use on highly parallel computers.  Because of the many matrix-vector products used per iteration with polynomial preconditioning, fewer iterations are needed and there can be a reduction in the orthogonalization expense.  With this comes a reduction in dot products and thus in communication and global synchronization for parallel linear solvers. There is potential to avoid even more communication expense by applying matrix-vector products using a Matrix Powers Kernel~\cite{AvoidComm,Ho10} or by combining polynomial preconditioning with pipelined methods. Furthermore, polynomial preconditioners are more suited for GPU computing than standard preconditioners such as incomplete factorization, which require sequential triangular solves.     


The GMRES polynomial is the polynomial implicit in the residual from one cycle of GMRES: $r = \pi(A)r_0$, where $r_0$ is the initial residual, $r$ is the residual computed at the end of the cycle, and $\pi$ is the GMRES polynomial.  The GMRES polynomial is a good choice, because it is relatively easy to compute and can be made stable.  The GMRES polynomial also adapts to the whole spectrum of $A$ instead of being based on estimates of the hull of the spectrum like a Chebyshev polynomial.  This could be advantageous for a matrix with a spectrum that has large gaps. 
The GMRES polynomial was previously investigated for a Richardson iteration by Nachtigal, Reichel, and Trefethen in~\cite{NaReTr} and by Joubert in~\cite{Jo} (he combines several GMRES polynomials to produce his polynomial for Richardson).  The GMRES polynomial was first used as a preconditioner for linear systems in~\cite{seed,PPG}.  The approach in~\cite{seed}, which is also the initial method in~\cite{PPG}, uses a power basis and normal equations to compute the coefficients of $p(\alpha)$.  
While this algorithm is very concise, it is only stable for low degrees; in some cases the calculation even returns NaN's for polynomial coefficients.  A more stable approach for $p$ is applied to eigenvalue problems in~\cite{Th06}.  However, in~\cite{PPG} it is shown that this and another more efficient attempt at a  stable approach can suffer from instability when an eigenvalue of $A$ is well-spaced from the others.  In this work, we give an approach that is more efficient than the more stable alternatives in~\cite{PPG} and is also much more stable.  Unlike in~\cite{PPG}, we focus on computing the $\phi(A)$ polynomial rather than the $p(A)$ polynomial, though computing $p$ is also needed.

Our new implementation of the GMRES polynomial preconditioner has two distinct components: We use one algorithm to apply $\phi(A)$ for (\ref{eqn:ppsys}) and a different but related approach to compute $p(A)$ in (\ref{eqn:ppsys2}).  We apply $\phi(A)$ using the formula $\phi(A) = I-\pi(A)$.  Here $\pi$ is the GMRES polynomial, and it is factored using its roots, which are the harmonic Ritz values~\cite{IE,PaPavdV,IEN, GoossensRitzValsConvGMRES,HRAM} generated via an initial GMRES run. This method of applying $\phi(A)$ is used in~\cite{PPArnoldi} to perform a spectral transformation for finding eigenvalues with the Arnoldi method.  The paper~\cite{PPArnoldi} also introduces the stability control method that we will use with linear systems: To improve the stability of the polynomial, we add extra copies of roots corresponding to eigenvalues that stand out in the spectrum.  Background information from \cite{PPArnoldi} will be presented in Section 2.  

Our work for linear equations deviates from the methods for eigenvalue problems in \cite{PPArnoldi} in several ways:  
\begin{enumerate}
\item We present an original algorithm for applying $p(A)$; details are in Section 3. (This application of $p(A)$ is not needed for eigenvalue problems.)  For the accuracy of the final computation of $x$ using equation~(\ref{eqn:ppsys2}), it is important to implement $p$ using an algorithm related to the one for $\phi$ (both use GMRES polynomial roots).  
\item A unique feature for linear systems is that the polynomial can be composed with a standard preconditioner such as incomplete LU factorization.  In that case, (\ref{eqn:ppsys}) and (\ref{eqn:ppsys2}) become 
\begin{align}
\phi(AM^{-1})y = b, \label{eqn:ppsyspre} \\  
x = M^{-1}p(AM^{-1})y,  \label{eqn:ppsys2pre}
\end{align}
where $M^{-1}$ is the standard preconditioner.  We give examples where polynomial preconditioning is used to accelerate an ILU preconditioned solve.  \item Section 3 demonstrates the root-adding stability method from \cite{PPArnoldi} applied to linear equations. Then we address stability considerations that were not presented in \cite{PPArnoldi}: We discuss the effect of adding roots on convergence.  A newly developed test detects unstable polynomials. Additionally, we show the importance of choosing a random starting vector (rather than the problem right-hand side) to generate $\phi$.   
\end{enumerate}

The remainder of our paper demonstrates the potential of polynomial preconditioning and establishes its differences from related methods.  Section 4 gives idealized estimates of the effectiveness of polynomial preconditioning for difficult problems and an example demonstrating this potential.  Section 5 has a comparison to FGMRES~\cite{Sa93} and a new variation of Polynomial Preconditioned (PP)-GMRES where we change the polynomial at every cycle.  Finally Section 6 has double polynomial preconditioning for additional reduction of dot products.

\section{Review}
\label{Sec:Review}
Here we review relevant methods from \cite{PPArnoldi} for applying the polynomial and for adding roots for stability.

GMRES works to solve $Ax=b$ by choosing an approximate solution $\hat{x}$ that minimizes the norm of the residual vector over the Krylov subspace $Span\{b, Ab, A^2b, \ldots,$ $A^{d-1}b\}$.  Thus, we can write $\hat{x} = p(A)b$ where the coefficients of $p$ correspond to the linear combination of Krylov vectors needed to form $\hat{x}$.  We can rewrite the residual vector as
\begin{equation*}
    r = b-A\hat{x} = (I-Ap(A))b = (I - \phi(A))b = \pi(A) b.
\end{equation*}
As GMRES builds a bigger subspace, the degree of $p(A)$ increases, the residual norm decreases, and ideally $p(A)$ becomes a better approximation to $A^{-1}$.  
Thus, we choose $\phi(A) = Ap(A)$ as our preconditioned operator.  The work \cite{PPArnoldi} runs one cycle of GMRES$(d)$ to find $\pi(A) = I - \phi(A)$ and then uses Arnoldi on the matrix polynomial $\pi(A)$ to compute eigenvalues and eigenvectors of $A$.  The matrix polynomial is implemented using 
\begin{equation}
    \pi(\alpha) = \prod_{i=1}^d\left(1-\frac{\alpha}{\theta_i}\right)
    \label{eqn:Qpoly}
\end{equation}
where the $\theta_i$'s are harmonic Ritz values, the roots of the $\pi$ polynomial.  The harmonic Ritz values are ordered with a modified Leja ordering \cite{BaHuRe} for numerical stability.  (If $A$ has complex entries, use a Leja ordering \cite{ReichelEvalPolyCoeffs} rather than a modified Leja ordering, and disregard the following information on avoiding complex arithmetic.)  

Since $A$ is real-valued, any complex harmonic Ritz values occur in conjugate pairs.  The modified Leja ordering sorts complex conjugates consecutively, allowing us to avoid complex arithmetic by combining conjugate pairs.  Suppose that $\theta_k = a + bi$ and $\theta_{k+1} = a - bi$.  Then 
\begin{align}
    \left( 1-\frac{\alpha}{\theta_k}\right)\left( 1-\frac{\alpha}{\theta_{k+1}}\right) 
    &= 1 + \frac{\alpha^2 -2\alpha a}{a^2+b^2}.
    \label{eqn:combineProd}
\end{align}
Thus, for $\pi(A)$ with a complex Ritz value, we can apply two factors of the polynomial together (see Algorithm \ref{Alg:APAtimesb}). 

\begin{algorithm}[]
    \caption{$\phi(A)$ times $v$}
    \Input  sparse matrix $A\in \R^{n\times n}$, $v \in \R^{n\times 1}$, $d$ harmonic Ritz values $\theta_i$. 
    \hspace*{\algorithmicindent}
    \begin{algorithmic}[1]
        \State $\pi = v$, $i=1$
            \While{$i<= d$}
                \If{$imag(\theta_i) == 0$}
                    \State $prod = A*\pi$
                    \State $\pi = \pi - (1/\theta_i)*prod$
                    \State $i=i+1$
                \Else
                    \State $a = real(\theta_i),\ b = imag(\theta_i)$
                    \State $prod = A*\pi$
                    \State $tmp = A*prod - 2*a*prod$
                    \State $\pi = \pi + (1/(a^2+b^2))*tmp$
                    \State $i=i+2$
                \EndIf
            \EndWhile
        \State $\phi = v - \pi$
        \State Return $\phi$.
    \end{algorithmic}
    \label{Alg:APAtimesb}
\end{algorithm}

Sometimes $\pi(\alpha)$ will have a very steep slope near one of the roots $\theta_k$.  This means that applying $\phi(A)$ to a vector can be ill-conditioned.  To resolve this problem in \cite{PPArnoldi}, $\pi(\alpha)$ is expanded to have extra copies of the term $(1-\alpha/\theta_k)$, i.e. to have a root of higher multiplicity at $\theta_k$.  This flattens the polynomial near $\theta_k$ and makes the preconditioner application more stable.  (See more details in Example 6 and Figure 3.2.)  To determine how many extra roots are needed at $\theta_k$, we compute a `product of other factors' or `$pof$' which estimates the slope of $\pi(\alpha)$ near $\theta_k$.  Experiments in \cite{PPArnoldi} suggested that one should add an extra root $\theta_k$ when $pof(k) > 10^4$ and another for every factor of $10^{14}$ beyond that.  Algorithm \ref{Alg:RootAdding} gives the procedure to automate adding roots for stability.  Subsections \ref{Sec:StabilityEx} and \ref{Sec:RootAdding} further discuss advantages and potential difficulties of adding roots.  Problems which benefit from added roots appear throughout the paper.


\begin{algorithm}
    \caption{Adding Roots to $\pi(\alpha)$ for Stability \cite{PPArnoldi}}
    \begin{enumerate}
 \item {\bf Setup:} Assume the $d$ roots ($\theta_1, \ldots, \theta_d$) of $\mrpoly$ have been computed and then sorted according to the modified Leja ordering~\cite[alg.~3.1]{BaHuRe}.  For very high degree polynomials, logs of products can be used to prevent overflow and underflow during the ordering.
 \item {\bf Compute $\pof(k)$:} For $k=1,\ldots,d$, compute $\pof(k) = \prod_{i \neq k} |1-\theta_k/\theta_i|$.  
 \item {\bf Add roots:} Compute least integer greater than $(\log_{10}(\pof(k)) - 4)/14$, for each $k$.  Add that number of $\theta_k$ values to the list of roots.  We add the first to the end of the list and if there are others, they are spaced into the interior of the current list, evenly between the occurrence of that root and the end of the list (keeping complex roots together).
\end{enumerate}
\label{Alg:RootAdding}
\end{algorithm}


\section{Polynomial Preconditioned GMRES using a GMRES polynomial}

\subsection{The method}
\label{Sec:PolyMethod}


As mentioned, the new approach for preconditioning linear systems lets $\phi(A) = I - \pi(A)$, where $\pi$ is the GMRES polynomial.  The $\pi$ polynomial is implemented using its harmonic Ritz value roots as in~\cite{PPArnoldi}.  This has about the same expense as the less-stable power basis approach in~\cite{PPG} but only about one-half 
the vector additions and dot products as in the cheapest attempt at a more stable method in~\cite{PPG}.  

Unlike with eigenvalue problems, an algorithm is also needed to apply $p(A)$ for the final step (\ref{eqn:ppsys2}).  
We give a method using harmonic Ritz values.  
Note in~\cite{Say88}, ways of finding roots of $p$ from those of $\pi$ are given, but only for special cases with three-term recurrences such as Chebyshev, so they do not apply here.  

To derive our approach, we start with the formula $\phi(A) = I-\pi(A)$ where $\pi(A)$ is written in factored form as in (\ref{eqn:Qpoly}).  Then we divide both sides by $\alpha$ and rewrite the polynomial $p(\alpha)$ as 
\begin{equation*}
    p(\alpha) = \frac{1}{\alpha} - \frac{1}{\alpha}\prod_{i=1}^d\left(1-\frac{\alpha}{\theta_i}\right).
\end{equation*} 
We split the last term of the product into two factors and distribute the $-1/\alpha$ to get
\begin{equation*}
    p(\alpha) = \frac{1}{\alpha} - \frac{1}{\alpha}\prod_{i=1}^{d-1}\left(1-\frac{\alpha}{\theta_i}\right) + \frac{1}{\theta_d}\prod_{i=1}^{d-1}\left(1-\frac{\alpha}{\theta_i}\right).
\end{equation*} 
Next we split the second term and distribute the $-1/\alpha$.  Continuing this process and canceling $1/\alpha$ terms at the end gives 
\begin{equation}
    p(\alpha) = \sum_{k=1}^d u_k \text{\ \  where\ \ }  u_k =  \frac{1}{\theta_k}\left( 1 - \frac{\alpha}{\theta_1}\right)\left( 1 - \frac{\alpha}{\theta_2}\right)\cdots\left( 1 - \frac{\alpha}{\theta_{k-1}}\right).
    \label{eqn:PDefn1}
\end{equation}
The algorithm for multiplying $p(A)$ by a vector alternates between building out the product for the next $u_k$ term and adding that term to the final sum.  


For real-valued matrices, we again avoid complex arithmetic by combining complex conjugates.  Suppose all $\theta_i$'s are real for $i <k$ and then $\theta_k = a+bi$ with $\theta_{k+1}=a-bi$.  Then the sum of the next two terms of $p(\alpha)$ is rewritten as follows: 
\begin{align*}
    u_k + u_{k+1} &=  \left( 1 - \frac{\alpha}{\theta_1}\right)\left( 1 - \frac{\alpha}{\theta_2}\right)\cdots\left( 1 - \frac{\alpha}{\theta_{k-1}}\right)\left( \frac{2a- \alpha}{a^2+b^2} \right). 
\end{align*}
Assuming $\theta_{k+1}$ is not the last root, we next need to form the product $u_{k+2}$. The last two terms of $u_{k+2}$ can be combined as in (\ref{eqn:combineProd}).  
Algorithm \ref{Alg:PAtimesbAltComplex} details the full process for applying $p(A)$ while avoiding complex arithmetic.  

\begin{algorithm}[h!]
    \caption{$p(A)$ times $v$}
    \Input  sparse matrix $A\in \R^{n\times n}$, $v \in \R^{n\times 1}$,  $d$ harmonic Ritz values $\theta_i$.  
    \hspace*{\algorithmicindent}
    \begin{algorithmic}[1]
        \State $prod = v$, $p = zeros(n,1)$, $i=1$ 
            \While{$i<=d-1$}
                \If{ $imag(\theta_i) == 0$}
                    \State $p = p +(1/\theta_i)* prod$
                    \State $prod = prod - (1/\theta_i)*A*prod$
                    \State $i=i+1$
                \Else
                    \State $a=real(\theta_i), \ b= imag(\theta_i)$
                    \State $tmp = 2*a*prod - A*prod$
                    \State $p = p + (1/(a^2 + b^2))*tmp$
                    \If{$i <= d-2$}
                        \State $prod = prod - (1/(a^2 + b^2))*A*tmp$
                    \EndIf
                    \State $i=i+2$
                \EndIf
            \EndWhile
            \If{$imag(\theta_d)==0$}
                \State $p = p + (1/\theta_d)*prod$
            \EndIf
        \State Return $p$.
    \end{algorithmic}
    \label{Alg:PAtimesbAltComplex}
\end{algorithm}

We let the vector operations or $vops$ be the total number of length-$n$ vector operations such as dot products, norms and $daxpys$ (a $daxpy$ is multiplying a vector by a scalar and adding to another vector).   
Note that applying $p(A)$ to a vector requires more $vops$ than applying $\phi(A)$. However, applying $p$ becomes less expensive as the number of complex harmonic Ritz values increases.  
If $r$ is the number of real harmonic Ritz values and $c$ is the number of non-real harmonic Ritz values, then the number of $daxpys$ needed to apply $p$ is $2r+(3/2)c-1$.  This gives one reason to compute $\phi(A)v$ using Algorithm 1 instead of using $A$ times $p(A)v$: if all the harmonic Ritz values are real, then applying $p$ directly requires almost twice as many $daxpys$ as applying $\phi$.  

Algorithm \ref{Alg:PPGMRESmain} summarizes the new polynomial preconditioned GMRES using the three previous algorithms.  To combine the polynomial with a standard preconditioner $M^{-1}$, simply use the matrix $AM^{-1}$ for the initial GMRES run and computation of harmonic Ritz values.  Then implement the second phase with (\ref{eqn:ppsyspre}) and (\ref{eqn:ppsys2pre}).


\begin{algorithm}[h!]
    \caption{GMRES with Polynomial Preconditioner of degree $\deg$}
\medskip
\begin{enumerate}
\item {\bf Construction of the polynomial preconditioner:} 
\begin{enumerate}
\item Run one cycle of GMRES($\deg$) using a random starting vector.
\item Find the harmonic Ritz values $\theta_1, \ldots, \theta_\deg$, which are the roots of the GMRES polynomial: With Arnoldi decomposition $AV_\deg = V_{\deg+1} H_{\deg+1,\deg}$, find the eigenvalues of $H_{\deg,\deg}^{} + h_{\deg+1,\deg}^2 @@f e_\deg^T$, where $f = H_{\deg,\deg}^{-*} e_\deg^{}$ with elementary coordinate vector $e_\deg=[0,\ldots,0,1]^T$.     
\item Order the GMRES roots and apply stability control as in Algorithm 2.
\end{enumerate}
\medskip
\item {\bf PP-GMRES:} Apply restarted GMRES to the matrix $\phi(A) = I-\Pi_{i=1}^{d}(I - A/\theta_i)$ to compute an approximate solution to the right-preconditioned system $\phi(A)y=b$, using Algorithm 1 for $\phi(A)$. To find $x$, compute $p(A)y$ using Algorithm 3.
\end{enumerate}
\label{Alg:PPGMRESmain}
\end{algorithm}




\subsection{Some experiments}
\label{Sec:E20Experiment}

This subsection has several examples of restarted GMRES with polynomial preconditioning.  The first three show that polynomial preconditioning can give a big improvement with difficult problems for which restarted GMRES converges slowly.  The next two give information to help discern cases where polynomial preconditioning can be effective.

We use GMRES(m), which restarts when the Krylov subspace reaches dimension $m$, both with and without the polynomial preconditioning.  All experiments have modified Gram-Schmidt orthogonalization with no reorthogonalization.  PP(d)-GMRES(m) refers to GMRES(m) with polynomial preconditioner of degree $d$.  When the degree is given with a plus sign, it means roots were added for stability, e.g. $150 + 2$ has original degree 150 and 2 added roots.  Unless stated otherwise, problem right-hand sides are generated random Normal(0,1) and then normed to one.  The initial guess is always $x_0 =\vec{0}$.  The experiments are run in Matlab on a Dell Optiplex desktop computer with I7-6700 processor.

{\it Example 1.}  We test the matrix E20r0100 from the Matrix Market collection.  It is nonsymmetric of size $n = 4241$ with an average of 31 non-zeros per row.  The corresponding linear equations are fairly difficult since the matrix is indefinite and has condition number $9.4*10^{6}$.  We run GMRES(50) with a residual norm convergence tolerance of $10^{-8}$.  
Table \ref{Tab:E20Ex} has results for polynomial preconditioning of $A$ (no standard preconditioning is used).  
The first row with $d = 1$ corresponds to no polynomial preconditioning, while other rows correspond to a $\phi$ polynomial of degree $d$ (so a $p$ polynomial of degree $d-1$).  The $mvps$ column indicates matrix-vector products, and $vops$ gives total length-$n$ vector operations.  The dot products included in the $vops$ count are given in a separate column.  The cycles column gives the number of restarts of GMRES$(50)$ plus $1$.  All results include the time and expense to create the preconditioner.   
Standard GMRES(50) does not converge.  Adding polynomial preconditioning improves results, though it does takes a high degree polynomial with $d = 150 + 2$ in order to get rapid convergence.  
The stability control with added roots keeps the residual norm from stalling before it reaches the requested tolerance. Without stability control, the problem with the polynomial of degree $150$ does barely reach the requested accuracy.
However, if extra roots are not added for degree 200, the residual norm only reaches $4.8*10^{-6}$.   
For this example, polynomial preconditioning is effective even with an indefinite spectrum.  However, we note that a very indefinite problem, especially with eigenvalues surrounding the origin, can be difficult for any Krylov method, including with polynomial preconditioning.  Another indefinite matrix is in Example 8, but we do not attempt to fully address indefinite problems in this paper.

\begin{table}

\caption{ Matrix E20r0100. Comparing no polynomial preconditioning to increasing degree polynomial preconditioning.  No standard preconditioning is used. }

\begin{center}
\begin{tabular}{|c|c|c|c|c|c|}  \hline\hline
degree  & cycles    & $mvps$       & $vops$        & dot products  & time          \\
d       &           & (thousands)   & (thousands)   & (thousands)   & (seconds)     \\  \hline \hline
1       & -         & -             & -             & -             & -             \\ \hline
25      & 215       & 268           & 860           & 285           & 33.6          \\ \hline
50      & 1231      & 3077          & 6468          & 1633          & 379           \\ \hline
100     & 204       & 1020          & 1591          & 275           & 109           \\ \hline
150 + 2 & 2         & 11.2          & 37.6          & 13.1          & 1.6          \\ \hline
200 + 4 & 2         & 10.6          & 54.0          & 21.6          & 1.9          \\ \hline \hline

\end{tabular}
\end{center}
\label{Tab:E20Ex}

\end{table}

\begin{figure}
\centering
\includegraphics[scale=.85]{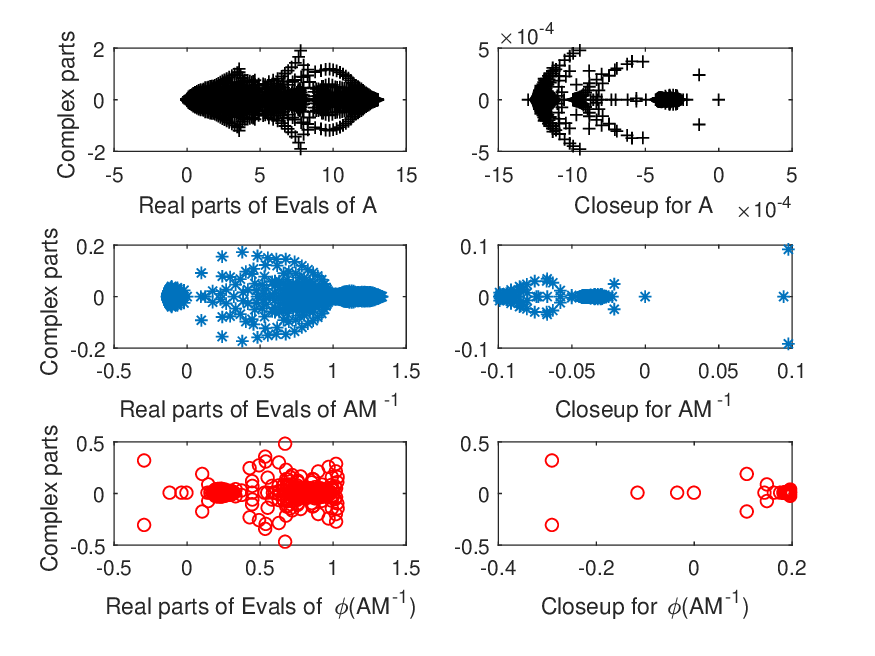}
\caption{Spectrum of matrix E20r0100 (top left) and a closeup of this spectrum near the origin (top right). The middle images show the spectrum after ILU$(0)$ preconditioning, and the bottom part has the spectrum after both ILU$(0)$ and polynomial preconditioning of degree $10$.}
\label{fig:E20Evals}
\end{figure}

We continue with matrix E20r0100 and add standard preconditioning.  Table \ref{Tab:E20ILU} has polynomials composed with several incomplete LU factorization preconditioners~\cite{MevdV,Sa03}.
The ILU$(0)$ factorization is from the shifted matrix $A + 0.01I$ (the ILU factorization of the original matrix $A$ fails due to zero on diagonal).  Then ILU(.01) uses the Matlab Crout ILU with fill-in tolerance of 0.01, also applied to $A + 0.01I$ (the factorization fails on $A$ due to zero pivot).  Finally, ILU(.001) is successfully applied to the original matrix $A$.
With ILU(0) preconditioning, GMRES(50) still does not converge until polynomial preconditioning is added, but now the low degree $d=10$ polynomial is effective. The combination of the two preconditionings makes the problem much easier for GMRES(50).  
Figure 3.1 shows the spectral transformation brought by the ILU$(0)$ and polynomial preconditionings.  The top of the figure has the spectrum of the matrix, along with a closeup showing that there are many negative eigenvalues.  The middle section has the spectrum with the ILU(0) preconditioning, and while it is vastly changed, there are still many negative eigenvalues.  Finally, the bottom part of the figure shows that if polynomial preconditioning of degree 10 is added to the ILU preconditioning, then the spectrum is much less indefinite and most of the eigenvalues are better separated from the origin.  In this example, polynomial preconditioning is effective even with a significantly complex spectrum.  
Also, even with standard preconditioning, polynomial preconditioning can be necessary to help the problem converge or can be used as an accelerator.

Next, we look at the ILU factorizations with fill-in.  
GMRES(50) finally does converge with ILU(.01), but adding polynomial preconditioning with $d=10$ makes it roughly four times faster in terms of matrix-vector products and time.  Vector operations are reduced by a factor of $34$.  ILU(.001) works well alone, and regular GMRES needs only 52 matrix-vector products.  This case shows that if standard preconditioning is very effective, the problem may become easy enough that polynomial preconditioning does not need to be added on.  When polynomial preconditioning is added, GMRES uses 91 matrix-vector products, but the solve time is about the same.  Vector operations, however, are still reduced by a factor of almost $4$. 

\begin{table}

\caption{  Matrix E20r0100 with several ILU preconditioners. }

\begin{center}
\begin{tabular}{|c||c|c|c||c|c|c|}  \hline\hline
 & \multicolumn{3}{|c||}{Regular GMRES} &  \multicolumn{3}{|c|}{With Poly, $d=10$}  \\ \hline
Type of ILU             & $mvps$    & $vops$        & time          & $mvps$        & $vops$        & time  \\
Preconditioner          &           & (thou's)      & (sec's)       &               & (thou's)      & (sec's)   \\ \hline \hline
ILU(0) of $A+.01 I$   & -         & -             & -             & 1001          & 6.44          & 0.65  \\ \hline
ILU(.01) of $A+.01 I$ & 1397      &   77.9        & 2.0           & 411           & 2.27          & 0.56  \\ \hline
ILU(.001) of $A$        & 52        &   2.71        & 0.52          & 91            & 0.71          & 0.58  \\ \hline
\hline

\end{tabular}
\end{center}
\label{Tab:E20ILU}
\end{table}

{\it Example 2.}
We use the fluid dynamics matrix Ill\_Stokes from SuiteSparse.  It is nonsymmetric with dimension $n=20{,}896$ and an average of 9.2 nonzeroes per row.  GMRES(50) is run to residual tolerance $10^{-8}$.  We test first polynomial preconditioning without regular preconditioning.  Then ILU$(0.001)$ is applied to the shifted matrix $A+0.001 I$ (this was one of the best ILU preconditioners we found).  Results are in Table \ref{Tab:IllStokes}.   With no regular preconditioning, polynomial preconditioning is essential as it reduces the computational time from 21.6 hours to 15 minutes (with polynomial of degree $100+20$).  With ILU preconditioning, it is still beneficial to add polynomial preconditioning.  The time is reduced from 211 seconds to about 14 seconds.  Matrix-vector products are reduced by a factor of about 7 and dot products by a factor of over 250.

\begin{table}

\caption{ Matrix Ill\_Stokes. Comparing no polynomial preconditioning to polynomial preconditioning both with and without standard preconditioning. }

\begin{center}
\begin{tabular}{|c|c|c|c|c|c|}  \hline\hline
degree d  & cycles  & $mvps$        & $vops$        & dot products  & time      \\  \hline \hline
\multicolumn{6}{|c|}{No Standard Preconditioning}    \\ \hline
1       & 485,042   & $2.43*10^7$   & $1.36*10^9$ & $6.43*10^8$ & 21.6 hours    \\ \hline
50 + 5  & 1072      & $2.95*10^{6}$ & $5.90*10^6$ & $1.42*10^6$ & 29.9 min's    \\ \hline
100 + 20 & 277      & $1.66*10^{6}$ & $2.43*10^6$ & $3.71*10^5$ & 15.2 min's    \\ \hline \hline
\multicolumn{6}{|c|}{With ILU Preconditioning}    \\ \hline
1       & 958       & 47,902        & $2.69*10^6$ & $1.27*10^6$ & 211 sec's             \\ \hline
50      & 3         & 7051          & 16,978    & 4799      & 13.6 sec's         \\ \hline
100 + 10 & 2        & 7691          & 21,273    & 6668      & 14.8 sec's          \\ \hline

\end{tabular}
\end{center}
\label{Tab:IllStokes}

\end{table}

{\it Example 3.}
We consider the 2-dimensional biharmonic partial differential equation $ - u_{xxxx} - u_{yyyy} +  u_{xxx} = f$ on the unit square, with $u=0$ on the boundaries.  The linear equations are challenging even without an extremely fine grid.  We use a uniform finite difference grid with both $\Delta x$ and $\Delta y$ of $\frac{1}{201}$, so that the matrix has size $n = 40{,}000$.  This matrix has 13 non-zeros in most rows.  
We stop when the shortcut residual norm formula reaches $10^{-10}$ (for this residual we use $||\beta e_1 - \bar H y||$, with the notation from~\cite{SaSc}).  The actual residual does not always reach this level due to the ill-conditioning of the problem.  Interestingly, while GMRES(50) reaches only $2.7*10^{-9}$ for the true residual norm, all the polynomial preconditioned tests reach $1.4*10^{-10}$ or better.
This matrix is very ill-conditioned with condition number $8.2*10^{7}$.  The polynomial preconditioning improves this substantially; for instance with $d=50$, the condition number of $\phi(A)$ is $6.1*10^{4}$.  
Table \ref{Tab:BihamrN40Th} has the results with some values of $d$.  The computational time is reduced by a factor of about 1100 going from no polynomial preconditioning to a polynomial of degree 403.  Looking at a breakdown of the expenses, the matrix-vector products go down by a factor of 156 for that polynomial (and even more for $d=819$).  The vector operations are lowest at $d=200$, with a reduction factor of 3700.  And the dot products go down even more, by a factor of almost 9000.  This is of particular interest for communication reduction.

\begin{table}

\caption{  Biharmonic matrix with $n=40{,}000$. }

\begin{center}
\begin{tabular}{|c|c|c|c|c|c|}  \hline\hline
degree  & cycles    & $mvps$       & $vops$        & dot products  & time          \\
d       &           & (thousands)   & (thousands)   & (thousands)   &               \\  \hline \hline
1       & 229,740   & 11,487        & 643,729       & 304,404       & 14.6 hours    \\ \hline
5       & 11,248    & 2,812         & 33,766        & 14,903        & 1.14 hours    \\ \hline
10      & 4,159     & 2,079         & 13,522        & 5,509         & 34.6 minutes  \\ \hline
25      & 742       & 927           & 2,969         & 983           & 10.4 minutes  \\ \hline
50      & 196       & 489           & 1,029         & 260           & 4.89 minutes  \\ \hline
100     & 47        & 235           & 374           & 137           & 2.29 minutes  \\ \hline
200     & 11        & 105           & 174           & 33.9          & 54.9 seconds  \\ \hline
400 + 3 & 4         & 73.7          & 245           & 85.1          & 47.9 seconds  \\ \hline
800 + 19& 1         & 41.7          & 688           & 322           & 1.33 minutes  \\ \hline
\hline

\end{tabular}
\end{center}
\label{Tab:BihamrN40Th}

\end{table}


It is reasonable to ask if there are certain classes of problems for which polynomial preconditioning for restarted GMRES will be particularly worthwhile.  We consider some possibilities, including problems of increasing difficultly and increasing non-normality.  

{\it Example 4.}  We use the same biharmonic matrix as in the previous problem, but add standard preconditioning.  We apply the incomplete factorization ILU(0) with the matrix first shifted in the positive direction by 0.5 times the identity matrix.  This shift is needed for the ILU to be effective (otherwise the factors are essentially singular).  Note that a more expensive factorization with fill-in was also not effective without the shift and did not improve results much with the shift.  This standard preconditioning makes the difficult problem easier, and so polynomial preconditioning is not needed as much, but still can be very helpful.  For the same matrix of size $n=40{,}000$, regular GMRES(50) now converges 368 times faster in terms of matrix-vector products with the ILU(0) preconditioning.  However with added polynomial preconditioning ($d=50$), it goes about another 13 times faster in $mvps$; see the next to last line of Table \ref{Tab:BiharmILU}.  Not shown in the table, dot products are reduced by a factor of 344 by adding polynomial preconditioning to the ILU preconditioning.

Next, we examine what happens as the difficulty of a problem increases.  
We continue with Table \ref{Tab:BiharmILU} and compare the results with and without polynomial preconditioning.  We start with a smaller matrix with $n=625$ and then increase the size and thus the difficulty.  For the matrix of size $n=625$, GMRES(50) needs fewer matrix-vector products and less time than PP(50)-GMRES(50), though vector operations are about the same.  As the difficulty increases, polynomial preconditioning passes up regular restarted GMRES, and by size $n=160{,}000$ it needs more than 20 times fewer matrix-vector products and almost two orders of magnitude less time.  
This points out that there is potentially a big advantage for polynomial preconditioning of difficult problems, while for easy problems, it may not be helpful.  

We now attempt to give more specific guidance for when polynomial preconditioning is worthwhile, at least for this situation.
Let $M=LU$ be from the incomplete factorization.
Then the condition numbers for $AM^{-1}$ of sizes $n=625$, $2500$, $10{,}000$, and $40{,}000$ are $9.4*10^{1}$, $1.1*10^{3}$, $1.6*10^{4}$, and $2.4*10^{5},$ respectively.  (The last one decreases to $1.9*10^2$ with $d=50$ polynomial preconditioning.) 
This gives an indication of how difficult the problem needs to be for the polynomial of degree $50$ in the table to be effective.  In terms of matrix-vector products, improvement occurs soon after $n=2500$, or soon after the condition number exceeds $10^{3}$.  In terms of computational time, polynomial preconditioning already gives a little speedup at this point.  For a lower degree polynomial, there can actually be improvement in matrix-vector products for the $n=2500$ case, as well.  With a polynomial of degree $5$, only 191 matrix-vector products are used.


\begin{table}

\caption{  Biharmonic matrix with ILU(0) preconditioning, changing matrix sizes. }

\begin{center}
\begin{tabular}{|c||c|c|c||c|c|c|}  \hline\hline
 & \multicolumn{3}{|c||}{Regular GMRES} &  \multicolumn{3}{|c|}{Polynomial prec with $d=50$}  \\ \hline
Size        & $mvps$    & $vops$        & time          & $mvps$        & $vops$        & time  \\
$n$     & (thousands)   & (thou's)      & (seconds)     & (thou's)      & (thou's)      & (sec's)   \\ \hline \hline
625         & 0.051     & 2.75          & 0.21          & 0.101         & 2.76          & 0.38  \\ \hline
2500        & 0.320     & 17.6          & 0.64          & 0.451         & 3.22          & 0.48  \\ \hline
10{,}000    & 2.81      & 160           & 7.6           & 0.751         & 3.69          & 0.86   \\ \hline
40{,}000    & 31.2      & 1780          & 208           & 2.30          & 7.25          & 8.5   \\ \hline
160{,}000   & 347       & 19{,}473      & $2.1*10^{4}$  & 16.1          & 36.1          & 233   \\ \hline
\hline

\end{tabular}
\end{center}
\label{Tab:BiharmILU}
\end{table}

{\it Example 5.}  
In this example, we look at the effect of increasing the non-normality of the matrix.  We also consider moving a complex spectrum closer to the left half of the complex plane.
We use a simple two-dimensional convection-diffusion equation $ - u_{xx} - u_{yy} + \alpha u_{x} + \beta u_{y} - \gamma u = f$ with zero boundary conditions on the unit square.   The matrix is generated with finite differences to size $n=40{,}000$.  No ILU preconditioning is used.

The symmetric problem with $\alpha=0$, $\beta=0$, and $\gamma=0$ is fairly difficult with condition number $1.6*10^4$.  Polynomial preconditioning is effective, reducing the time by about a factor of 12; see the first line in Table \ref{Tab:ConvDiff}.  As convection coefficients $\alpha,$ and $ \beta$ increase to $25$, then to 100 and 400, (see next lines of the table), the problem becomes easier: The eigenvalues become complex, but move away from the origin.  When $\alpha = \beta = 25$, for instance, the condition number goes down to $5.5*10^3$.  Even though condition numbers go down further for the larger $\alpha$'s and $\beta$'s, the results are about the same (apparently the movement of the spectrum further into the complex plane balances the movement away from the origin).  Polynomial preconditioning is not as effective for these non-normal but easier problems, but still reduces the solve time by more than a factor of four.  

In the last five lines of the table, we fix the convection coefficients at 400; here the spectrum is roughly circular in shape.  When $\gamma=0$, the smallest eigenvalues are $0.15 \pm 0.52 i$, which sit well away from the imaginary axis (relative to the largest eigenvalue, which is near 8).  To make the problem more difficult, we increase the value of $\gamma$ so that the spectrum approaches the origin and almost encircles it above and below.  For this non-normal matrix, regular GMRES(50) has difficulty well before the spectrum gets to the origin, because the pseudospectrum~\cite{TrEm} surrounds the origin first.  However, PP(50)-GMRES(50) is less affected by this.  Without polynomial preconditioning, the method does not converge when $\gamma$ reaches 1070.  On the other hand, with polynomial preconditioning there is convergence even at $\gamma = 10{,}000$.  As $\gamma$ further increases, convergence becomes slower and stalls by $\gamma = 12{,}500$.  
We conclude that polynomial preconditioning can assist both normal and quite non-normal matrices.

\begin{table}

\caption{ Convection-diffusion equation matrix with $n=40{,}000$ and with various coefficients. }

\begin{center}
\begin{tabular}{|c||c|c|c||c|c|c|}  \hline\hline
 & \multicolumn{3}{|c||}{Regular GMRES} &  \multicolumn{3}{|c|}{Polynomial prec with $d=50$}  \\ \hline
Coefficients    & $mvps$    & $vops$        & time          & $mvps$        & $vops$        & time  \\
$\alpha$, $\beta$, $\gamma$ & (thou's) & (thou's)  & (sec's)    & (thou's)      & (thou's)      & (sec's)   \\ \hline \hline
0, 0, 0         & 3050      & 171           & 12.8          & 1051          & 4.16          & 0.95  \\ \hline
25, 25, 0       & 950       & 53.0          & 4.3           & 951           & 3.98          & 0.96  \\ \hline
100, 100, 0     & 1018      & 56.4          & 4.6           & 951           & 3.98          & 0.96  \\ \hline
400, 400, 0     & 942       & 52.2          & 4.3           & 951           & 3.98          & 0.97  \\ \hline
400, 400, 1000  & 1934      & 83.1          & 8.7           & 1001          & 4.07          & 0.96  \\ \hline
400, 400, 1060  & 3462      & 194           & 14.4          & 1051          & 4.16          & 0.97  \\ \hline
400, 400, 1070  & -         & -             & -             & 1051          & 4.16          & 0.96  \\ \hline
400, 400, {10,000}  & -     & -             & -             & 4901          & 12.7          & 2.2  \\ \hline
\hline

\end{tabular}
\end{center}
\label{Tab:ConvDiff}
\end{table}

\subsection{Stability}
\label{Sec:StabilityEx}
Stability control is typically needed for high degree polynomials when there is an eigenvalue that stands out from the rest of the spectrum.  In this situation, the polynomial will have large slope at the eigenvalue.  This slope gives ill-conditioning and causes significant numerical error and a lack of convergence.  However, extra roots can control the steep slope.  

{\it Example 6.}  We look at stability control with the matrix $1138\_$bus from Matrix Market.  It is symmetric of size $n=1138$ with an average of 3.6 non-zeros per row.  No standard preconditioning is used.  Results in Table \ref{Tab:StabEx} have stability control on the left side and no added roots on the right side.  We first observe that this is a tough matrix, and as a result, polynomial preconditioning has room to be very effective.  Solve time is reduced by a factor of 400 for no preconditioning versus a polynomial with either degree $50 + 15 = 65$ or degree $75+61=136$.  It is interesting that the stability procedure in Algorithm 2 adds such a large number of extra roots.  This example needs many added roots because the larger eigenvalues are in bunches well-separated from the others.  With stability control, the results stay accurate even for the degree $75+61$ polynomial (accuracy does degrade after that point, as will be mentioned in Example 9).  We see in the right-most table column that without stability control, accuracy decreases as the degree increases and by $d=35$ there is essentially no accuracy (beyond that degree, the iterations do not converge, even in the shortcut residual.)   

\begin{figure}
\hspace{-.30 in}
\includegraphics[scale=0.95]{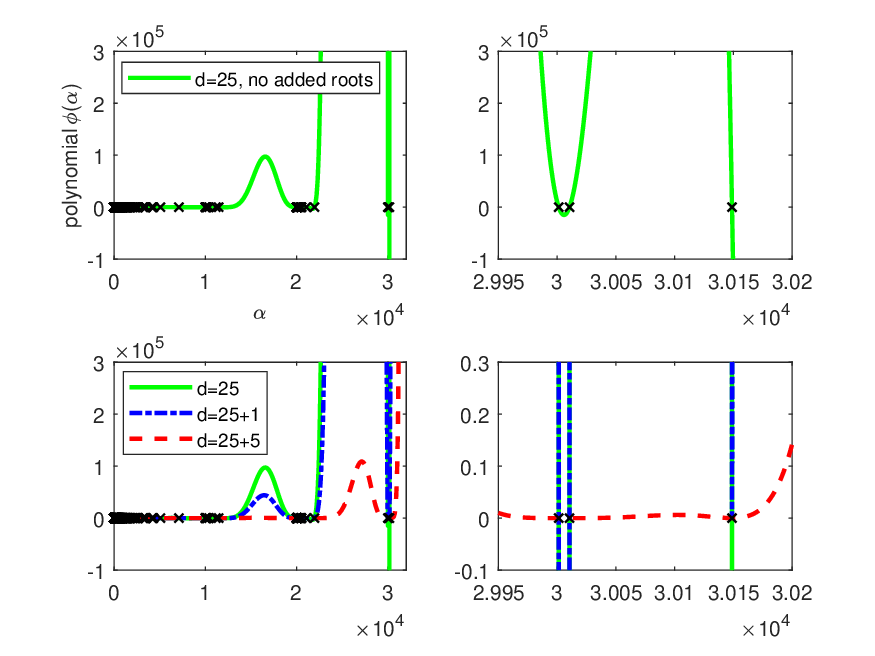}
\caption{The matrix is $1138\_$bus.  Graphs of the degree 25 GMRES polynomial (solid line), the polynomial with one root added at 30{,}149 (dash-dot line) and with five added roots (dashed line).  Eigenvalues are marked with x's.}
\label{Fig:StabEx}
\end{figure}

We now look at how the polynomial changes with added roots and why this makes a better conditioned polynomial and, thus, a more stable method.  The upper left portion of Figure~\ref{Fig:StabEx} has a graph of the GMRES polynomial of degree 25 and has the eigenvalues of $A$ marked with x's.  At first glance, it appears that there is a single eigenvalue at 3 that the polynomial goes down through steeply.  However, the closeup in the upper right shows there are three eigenvalues near 3.  The polynomial is very steep at the largest one and somewhat steep at the other two (more so than appears in the figure due to scaling).  The steep slopes cause the polynomial to be ill-conditioned.  In the lower left portion of the figure, the dash-dot (blue) line shows the degree $25+1=26$ polynomial with one root added at the largest eigenvalue.  It has slope zero at the large eigenvalue due to the double root, but slopes are still steep at the two nearby eigenvalues.  The polynomial of degree $25+5=30$ formed by the stability control, shown with dashed (red) line, has roots added at all three of the eigenvalues near 3 and two roots added near 2.  It has much less steep slopes.  The scale on the y-axis of the lower right graph is six orders of magnitude smaller than that of the other graphs; this shows how much milder the slopes are for the degree 30 polynomial.

\begin{table}

\caption{Matrix $1138\_$bus with and without stability control.}

\begin{center}
\begin{tabular}{|c|c|c|c|c||c|c|c|}  \hline\hline
& \multicolumn{4}{|c||}{With Stability Control} &  \multicolumn{3}{|c|}{Without Stab. Control}      \\ \hline
degree      &  added        & $mvps$    & time      & res. norm         & $mvps$    & time      & res. norm         \\ 
            &  roots        & (thou's)  & (sec's)   &                   & (thou's)  & (sec's)   &                   \\ \hline 
1           &               & 696       & 159       & $1.2*10^{-10}$    & 696       & 159       & $1.2*10^{-10}$    \\ \hline
15          & + 1           & 66.2      & 2.3       & $3.8*10^{-12}$    & 59.8      & 2.1       & $6.3*10^{-12}$    \\ \hline
25          & + 5           & 12.3      & 0.50      & $5.2*10^{-12}$    & 106       & 3.2       & $1.4*10^{-6}$     \\ \hline 
30          & + 7           & 15.3      & 0.53      & $5.6*10^{-12}$    & 65.7      & 1.9       & $5.4*10^{-3}$     \\ \hline 
35          & + 8           & 16.8      & 0.57      & $6.9*10^{-12}$    & 14.6      & 0.53      & $3.3*10^{+1}$     \\ \hline 
50          & + 15          & 9.87      & 0.39      & $8.1*10^{-12}$    & -         & -         & -                 \\ \hline
75          & + 61          & 7.96      & 0.36      & $1.1*10^{-11}$    & -         & -         & -                 \\ \hline 
\hline

\end{tabular}
\end{center}
\label{Tab:StabEx}

\end{table}

\subsection{Effect of adding roots}
\label{Sec:RootAdding}

Here we examine the effect of adding roots on the convergence rate.  The added roots require extra matrix-vector products per outer GMRES iteration, but if located to the outside of the spectrum, they typically make the polynomial better.  Of course, some slowing of convergence is acceptable if the added roots make the answer more accurate.  We start with the added roots in Example 6.

{\it Example 6 (cont.)}  The costs ($mvps$ and time) for matrix $1138\_$bus are compared with and without the stability control.  The data in Table \ref{Tab:StabEx} is for achieving residual tolerance of $10^{-10}$ in the shortcut residual formula, not the full accuracy that is eventually attained.  For original degrees 15 and 35, the method without added roots is slightly more efficient, so this shows that adding roots can reduce the effectiveness.  However, for original degrees 25 and 30, the controlled method is much better.  And even with 61 roots added to the degree 75 polynomial, it still gives the best time of all tests.

The following simple theorem shows that adding an extra root in the larger half of the spectrum does not hurt a Richardson iteration with the GMRES polynomial.  

\begin{theorem}
Assume $A$ is symmetric, positive definite with largest eigenvalue $\lambda_n$.
The Richardson iteration with the GMRES polynomial applied to $Ax=b$ computes $\pi(A)b$, where $\pi(\alpha) = \prod_i^n (1 - \frac{\alpha}{\theta_i})$.  Adding a term $(1 - \frac{\alpha}{\theta_k})$ to the product, where $\theta_k$ is a root of the GMRES polynomial and $\theta_k \ge \frac{\lambda_n}{2}$, does not cause the residual norm to increase.  In other words, $\|(I - \frac{1}{\theta_k}A) \pi(A)b\| \le \|\pi(A)b\|.$
\end{theorem}

{\em Proof}.  
It is given that $\theta_k \ge \frac{\lambda_n}{2}$.  (We also know that $\theta_k \le \lambda_n$, because $\theta_k$ is a harmonic Ritz value and $A$ is SPD, but this is not needed).  The line $y=(1 - \frac{\alpha}{\theta_k})$ satisfies $|y| \le 1$ on $[0, \lambda_n]$ and thus on the spectrum of $A$, and this gives the result.  
\endproof

We do not have a corresponding result for polynomial preconditioned GMRES, but the outer GMRES for PP-GMRES does use this Richardson polynomial at each iteration.

Next we look at how there can be a problem with adding roots: In rare cases, it can create a preconditioned problem that is worse for GMRES.  We have not seen this happen in definite application matrices that we have tested, because this requires large gaps among the smaller eigenvalues (see Example 8 for an indefinite example).  For the SPD case, this cannot happen when only the smallest eigenvalue is isolated, because a root is never added for the smallest harmonic Ritz value: If $\theta_1$ is the smallest in magnitude harmonic Ritz value and $\theta_i$ is any other, then $|1 - \frac{\theta_1}{\theta_i}| < 1$. 
This means that in Algorithm 2, $pof(1) < 1$ and will not trigger an extra root.
However, it is possible to construct an SPD example with small (but not the smallest) isolated eigenvalues and a high degree polynomial where adding roots has a bad effect.

\begin{figure}
\centering
\includegraphics[scale=0.80]{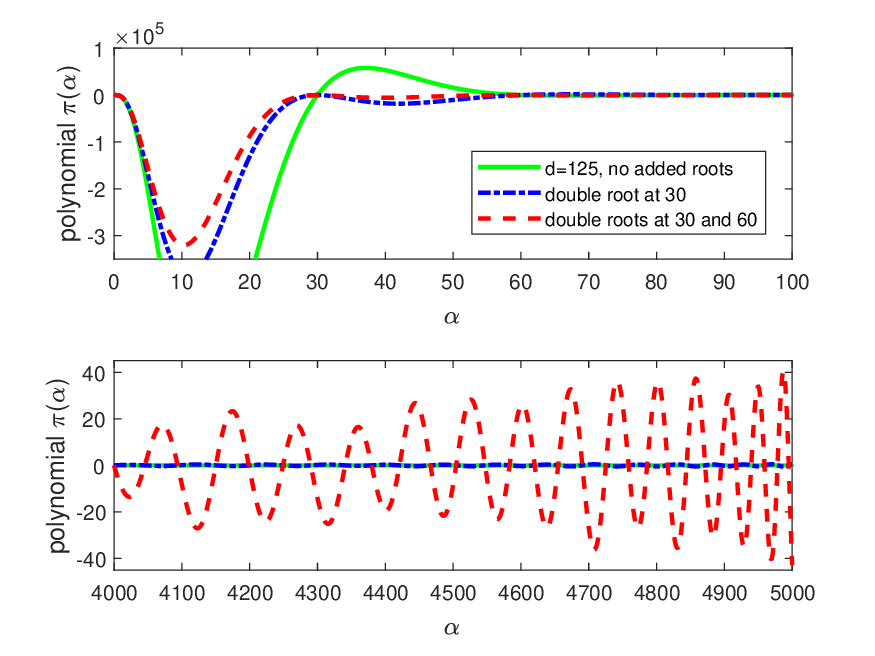}
\caption{For the diagonal matrix of Example 7, plots show the degree 125 GMRES polynomial (solid line), the polynomial with one root added at 30, (dash-dot line) and then with two roots added at 30 and 60 (dashed line). The top shows the polynomials in the small part of the spectrum and the bottom has the increased magnitude of oscillations in the large part of the spectrum when roots are added.}
\label{Fig:BadRootsAdd}
\end{figure}

{\it Example 7.}
Let $A$ be the diagonal matrix with entries $0.01, 0.02, \ldots 0.09, 0.1, 30,$ $60,$ $100, 101, \ldots 5087$ and $n= 5000$.
We run PP(125)-GMRES(50) and because of the well-separated eigenvalues, the $d=125$ polynomial has two roots added (one at $30$ and one at $60$).  These added roots in the lower part of the spectrum cause the polynomial to become large in the higher part of the spectrum.  The top of Figure~\ref{Fig:BadRootsAdd} shows the original GMRES polynomial in solid (green), and it has steep slope at the eigenvalue at $30$.  Adding one root at $30$ (blue dash-dot) tames this slope, and adding roots at both $30$ and $60$ (red dashed) makes it even flatter in that region.  The lower half of the figure shows how the GMRES polynomial $\pi$ with two added roots has large oscillations in the high part of the spectrum.  Before adding roots, the polynomial oscillates between approximately $-0.003$ and $0.003$ in the larger part of the spectrum.  With just the one root added at $30$, the polynomial goes between $-0.5$ and $0.5$, and the method still converges, but four times slower than without the added root.  With the two added roots, the oscillations are much bigger, and this turns the spectrum of $\phi(A)$ badly indefinite and the convergence is extremely slow.  
A possible solution for this problem is to use a lower degree polynomial.  Also, the added roots are not actually needed to get residual norm below $10^{-10}$, so the $pof$ cutoff in Algorithm 2 can be adjusted if higher accuracy is not required (the $pof$’s for the degree 125 polynomial are 5.79e+05 and 3.12e+04 for 30 and 60, respectively).

\subsection{Polynomial stability check}
\label{Sec:StabilityTest}

Even with stability control, there is a possibility that a high degree polynomial will be unstable.  Here we give a test that can suggest whether a particular polynomial will be stable.  This stability check can be applied once the polynomial has been determined, before the PP-GMRES linear solve.  Then the degree can be lowered if instability is predicted.  

We compute a residual norm for a very rough approximate solution $\hat x = p(A)b$ in two ways and compare them.  First $r_1 = b - A \hat x = b - A p(A) b$, where $p$ is implemented with Algorithm 3.  Then $r_2 = b - \phi(A)b = \pi(A)b$, using the factored form of $\pi$ in Algorithm 1.  Then the stability check is $StCh = \|r_1 - r_2\|$.  $StCh$ gives an estimate of the limit of the residual norm convergence to be expected in the PP-GMRES phase.  In our testing, the actual error is generally within an order of magnitude or two of $StCh$.

{\it Example 8.} We use the matrix OLM1000 from Matrix Market.  It has size $n=1000$.  It has some complex eigenvalues and is a little indefinite; all but $10$ eigenvalues have negative real parts. 
We look at polynomial preconditioning after first applying ILU(0) preconditioning.  The preconditioned spectrum is still indefinite, so extra roots added on the left side of the spectrum can increase the slopes of the polynomial on the right side.  Thus, the stability control can actually cause instability.  More research is needed for indefinite problems, 
but for now we merely show that this difficulty can be detected before the linear solve.  Table \ref{Tab:StabCh} has some choices of polynomials.  The first two have $d=10$ and 12 without the additional roots indicated by Algorithm 2, so they are ill-conditioned.  Then the last three polynomials have roots added, but still can have instability as the degree increases.  These tests all show reasonable correspondence between the value of $StCh$ and the residual norm at the end of the linear solve.

\begin{table}

\caption{Stability check versus minimum residual norm attained for OLM1000 with ILU(0) preconditioning.  First two polynomials are without added roots and the others do have added roots.}

\begin{center}
\begin{tabular}{|c|c|c|c|c|c|}  \hline\hline
degree      & 10            & 12            & 25                & 30            & 35            \\
            & (no added)    & (no added)    & + 6               & + 8           & + 15          \\ \hline \hline
StCh        & $1.2*10^{-5}$ & $5.4*10^{-3}$ & $4.3*10^{-11}$    & $5.0*10^{-8}$ & $1.9*10^{-1}$ \\ \hline
Res. Norm   & $3.6*10^{-5}$ & $2.6*10^{-2}$ & $1.4*10^{-10}$    & $3.0*10^{-8}$ & $6.0*10^{-1}$ \\ \hline
\hline

\end{tabular}
\end{center}
\label{Tab:StabCh}

\end{table}

{\it Example 9.}  We return to the matrix $1138\_$bus from Example 6.  With polynomial of original degree 100, there are 115 roots added for stability, however there is significant loss of accuracy.  The actual residual norm reaches minimum value of $2.1*10^{-5}$ and varies slightly at each iteration, so it can be higher depending on which iteration one stops at.  The reason that there are inaccurate results is probably due to the ordering of roots given in Algorithm 2.  The scalar polynomial with this ordering is not accurate when evaluated near the largest eigenvalue.  Perhaps this could be fixed with a better ordering of the roots, and this should be studied in the future.
Table \ref{Tab:StabChbus} has both the minimum residual norm and the value of StCh for different degree polynomials.  StCh again gives an indication of when there are problems with the polynomial.  Note that in Examples 7 and 8, there are three different sources of error, and the StCh test works for all of them.

\begin{table}

\caption{Stability check versus minimum residual norm attained for $1138\_$bus. }

\begin{center}
\begin{tabular}{|c|c|c|c|c|c|}  \hline\hline
degree      & 75        & 90            & 100               & 110           & 120           \\
            & + 61          & + 81          & + 115             & + 174         & + 236         \\ \hline \hline
StCh        & $2.3*10^{-13}$& $1.0*10^{-9}$ & $1.9*10^{-6}$     & $2.4*10^{-5}$ & $2.3*10^{+1}$ \\ \hline
Res. Norm   & $1.1*10^{-11}$& $9.0*10^{-9}$ & $2.1*10^{-5}$     & $8.1*10^{-4}$ & $7.4*10^{+1}$ \\ \hline
\hline

\end{tabular}
\end{center}
\label{Tab:StabChbus}

\end{table}

\subsection{The starting vector for the polynomial}
\label{Sec:StartVec}
When generating the polynomial preconditioner, it is typically best to run GMRES$(d)$ with a random right-hand side rather than using the problem right-hand side.  The following experiment demonstrates how polynomials generated using the problem right-hand side might ignore certain eigenvalues and give bad preconditioners.

\begin{figure}
    \centering
    \includegraphics[scale=.75]{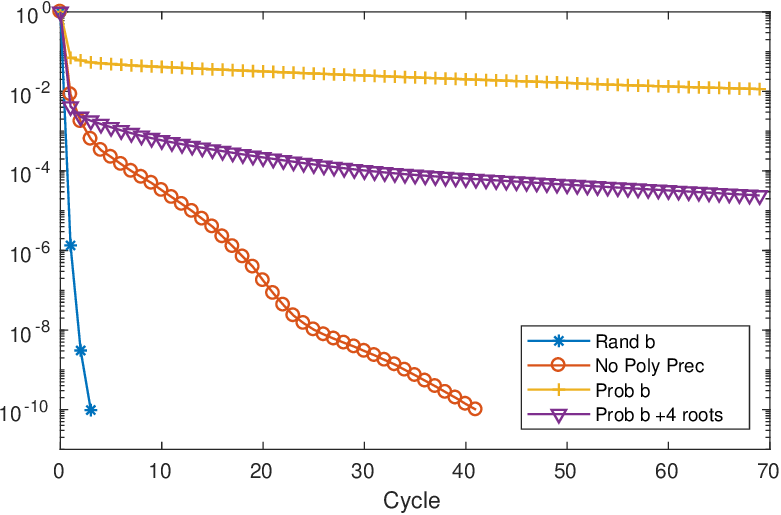}
    \caption{Relative Residual convergence vs number of cycles for GMRES(50) on the Memplus matrix.  Circles indicate no polynomial preconditioning.  The preconditioned problems are indicated by: asterisks for the $b_{rand}$ polynomial (degree $d=15$), crosses for the $b_{prob}$ polynomial ($d=15$), and triangles for the polynomial with added roots ($d = 15 + 4$). }
    \label{fig:memplusResNorms}
\end{figure}

{\it Example 10.}  We consider the electronic circuit matrix Memplus and its corresponding right-hand side available on Matrix Market.  The matrix $A$ is of size $n=17{,}758$. We let $b_{prob}$ denote the problem right-hand side and $b_{rand}$ denote a vector generated from a random Normal$(0,1)$ distribution.  
We generate two polynomial preconditioners of degree $d=15$.  The first is created by running GMRES(d) with $b_{rand}$ as a starting vector and the second by using $b_{prob}$ as a starting vector.  For the second polynomial, Algorithm 2 indicates that four additional roots are needed for stability, and this gives a third polynomial with $d=15+4$.  Figure \ref{fig:memplusResNorms} shows GMRES$(50)$   residual norm convergence while solving $Ax=b_{prob}$.  Tests are without polynomial preconditioning and with the three polynomials.  While the problem does converge in $41$ cycles without a preconditioner, the $b_{rand}$ polynomial gives a $90\%$ decrease in $daxpys$ and a $94\%$ decrease in dot products.  This comes at a price of a slight increase in $mvps$ with the polynomial.  However, the polynomial generated with $b_{prob}$ stalls GMRES convergence and makes the problem much worse than with no preconditioning.  Since instability is not the dominant problem with this ineffective polynomial, added roots give little improvement.   

\begin{figure}
    \begin{subfigure}[b]{.49\linewidth}
        \centering 
        \includegraphics[scale=.63]{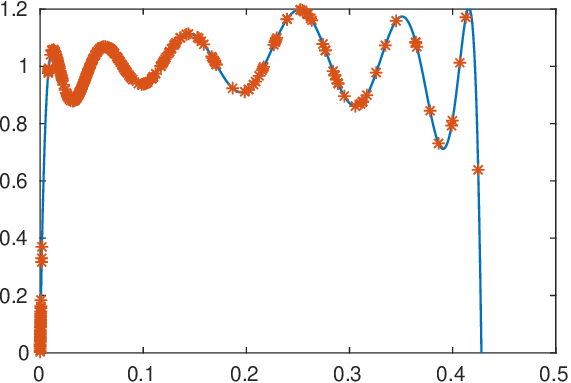}
        \caption{Polynomial generated with $b_{rand}$}
        \label{fig:memplusGoodPoly}
    \end{subfigure}
    \begin{subfigure}[b]{.49\linewidth}
        \centering
        \includegraphics[scale=.63]{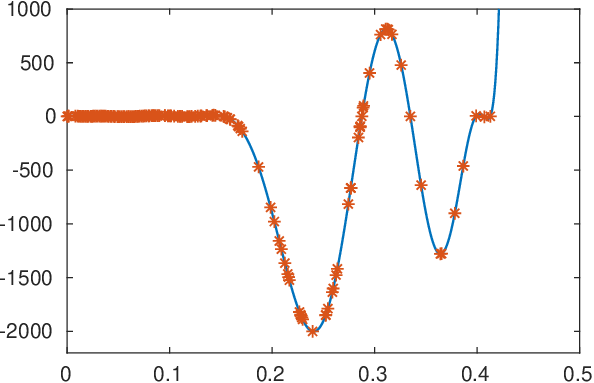}
        \caption{Polynomial generated with $b_{prob}$}
        \label{fig:memplusBadPoly}
    \end{subfigure}
    \caption{The $\phi$ polynomials of degree $15$ plotted over the real axis on $[0,0.5]$.  Stars indicate eigenvalues of the Memplus matrix (horizontal axis) mapped to the eigenvalues of the preconditioned matrix (vertical axis). Observe the large difference in scaling between the two vertical axes.}
    \label{fig:MembplusGoodvsBadPoly}
\end{figure}

This inconsistent preconditioner behavior can be explained by considering how the polynomials remap the eigenvalues of $A$.  All eigenvalues of $A$ lie in the right half of the complex plane. The two non-real eigenvalues have very small magnitude and are minimally affected by the polynomial preconditioners.  Of the real eigenvalues, four lie near $1.5$, and the remainder have magnitude less than or equal to $0.5$. Figure \ref{fig:MembplusGoodvsBadPoly} shows the effect of the first two polynomial preconditioners on this subset of real eigenvalues (note that though the spectrum is slightly complex, the polynomial is only graphed on the real axis).  While the $b_{rand}$ polynomial effectively maps most of the small eigenvalues to near $1$, the $b_{prob}$ polynomial creates a more difficult spectrum by making the problem highly indefinite.  Outside of this figure, the $b_{rand}$ polynomial moves the eigenvalues of magnitude $1.5$ closer to $1$, but the $b_{prob}$ polynomial effectively ignores those eigenvalues, mapping them to near $10^6$.  

Upon further examination of the vector $b_{prob}$, it appears that $b_{prob}$ has components of significant magnitude in the direction of only a handful of the eigenvectors of $A$. The components of $b_{rand}$ are more evenly distributed among the eigenvectors.  
Experiments in generating polynomials with other vectors yielded comparable results: A uniformly distributed random vector worked almost as well as the $b_{rand}$ vector, but a vector of all ones again yielded stalled convergence.  
This phenomenon has been observed with a number of other matrix problems.  Thus, we recommend always using a random vector to generate the polynomial preconditioner.  

\section{Potential of Polynomial Preconditioning}

In this section, we look at potential effectiveness of polynomial preconditioned GMRES for difficult problems.  
We use Chebyshev estimates to develop a theoretical estimate for how effective polynomial preconditioning can be for improving restarted GMRES.  There is dramatic reduction in the number of matrix-vector products under idealized circumstances.

We assume that all of the polynomials in polynomial preconditioned GMRES can be approximated with Chebyshev polynomials.  This includes both the polynomial for the preconditioning and the polynomials that underlie the GMRES method.  We assume all the eigenvalues of $A$ are
real and positive and lie between $a$ and $b$ with $0 < a < b $.  We also assume that the linear equations problem is very difficult, so $b$ is much larger than $a$.  As a result of this assumption, it is possible to use the approximation \begin{equation}
T_m(1+\delta) \doteq 1 + m^2 \delta, 
\label{eqn:Chebappr}
\end{equation}
where $T_m$ is the standard Chebyshev polynomial of the first kind of degree $m$ and $\delta$ is very small.
For one cycle of GMRES(m), the GMRES polynomial is assumed to be approximately the Chebyshev polynomial shifted and scaled so that it is one at the origin and small and equal oscillatory over the interval from $a$ to $b$.  Then the maximum value of the Chebyshev polynomial over the interval $[a, b]$ is $\frac{1}{T_m(1+2 a/(b - a))}$.  This quantity then gives approximately how much one cycle of GMRES(m) improves the residual norm.  With the approximation in (\ref{eqn:Chebappr}), we have that the residual norm is improved by approximately
\begin{equation*}
\frac{1}{1+\frac{2 m^2 a}{b - a}} \doteq 1-\frac{2 m^2 a}{b - a} \doteq 1-\frac{2 m^2 a}{b} .
\end{equation*}

We next compare the improvement in residual norm for $d$ cycles of GMRES(m) with the improvement for one cycle of polynomial preconditioned GMRES(m) with polynomial of degree $d$.  The improvement factor for the $d$ cycles of GMRES(m) is approximately
\begin{equation}
\Big(1-\frac{2 m^2 a}{b} \Big)^d \doteq 1-\frac{2 d m^2 a}{b} .
\label{eqn:ChebyEst1}
\end{equation}
We view the one cycle of polynomial preconditioned GMRES as being a composition of two polynomials with the preconditioner polynomial from GMRES(d) on the inside and a GMRES(m) poly on the outside.  This can be modeled with a composition of shifted and scaled Chebyshev polynomials, giving residual improvement of
\begin{equation}
\frac{1}{T_m(T_d(1+\frac{2a}{b - a}))} \doteq \frac{1}{T_m(1+\frac{2 d^2 a}{b - a})} \doteq \frac{1}{1+\frac{2 d^2 m^2 a}{b - a}} \doteq 1-\frac{2 d^2 m^2 a}{b}.
\label{eqn:ChebyEst2}
\end{equation}
Comparing (\ref{eqn:ChebyEst1}) and (\ref{eqn:ChebyEst2}), we conclude that polynomial preconditioned GMRES converges approximately $d$ times faster.  This is summarized in the following estimate. 
\begin{estimate}
For $A$ with real spectrum between $a$ and $b$ and with $0 < a \ll b $,  PP(d)-GMRES(m) converges approximately $d$ times faster than GMRES(m) in terms of matrix-vector products.
\end{estimate}

So where this estimate holds, we expect that if the degree of the polynomial is doubled, then the number of matrix-vector products is cut in half.  Note that orthogonalization costs are reduced even more dramatically.

Since this result uses two kinds of approximations (GMRES polynomials are approximated by Chebyshev polynomials and values of Chebyshev polynomials are approximated with an asymptotic result), it is natural to ask whether such a reduction in matrix-vector products can happen in a computation.  The next example tests this for a simple but difficult problem.

{\it Example 11.}  We let the matrix be diagonal with entries $\frac{1^2}{n}, \frac{2^2}{n}, \frac{3^2}{n}, \ldots, \frac{n^2}{n}$, where $n = 20{,}000$.  The residual tolerance is $10^{-10}$.  GMRES(50) is used both with and without polynomial preconditioning.   The degree of the polynomial is doubled for each test and we look at how the matrix-vector products are reduced.  The last two tests have roots added for stability.  Table \ref{Tab:TheoryScaling} shows the results, and while the matrix-vector products are not quite cut in half for each subsequent test, they do come close.  The last experiment with $d = 1048$ does reduce the matrix-vector products by a factor of 1360 compared to $d=1$, no polynomial preconditioning.  

Although they are not included in Estimate 4.1, it is interesting to note the reductions in solve time and other operations in this example. The $vops$ go down by a factor of over 8000 from $d=1$ to $d=516$.  In this ideal situation with both a difficult problem and cheap matrix-vector products so that most of the expense is in the orthogonalization, the reduction in computational time is remarkable for these degrees.  The time goes down from 3.85 days to 24.8 seconds, a reduction by a factor of about $13{,}400$.  The time then goes up for $d = 1048$ due to the initial cost of computing the polynomial using one cycle of GMRES(1024). For high-performance computing, dot products are expected to be a bottleneck, and they are reduced by a factor of over 20,000 from $d=1$ to $d=256$.

\begin{table}

\caption{ Example to demonstrate reduction in matrix-vector products that is roughly proportional to the degree of the polynomial.  $A$ has $n=20{,}000$ and is diagonal with entries $\frac{i^2}{n}$, $i=1 \ldots n$. }

\begin{center}
\begin{tabular}{|c|c|c|c|c|c|}  \hline\hline
degree  & cycles    & $mvps$       & $vops$        & dot products  & time          \\
d       &           & (thousands)   & (thousands)   & (thousands)   &               \\  \hline \hline
1       & 1,395,850 & 69,792        & 3,911,170     & 1,849,500     & 3.85 days     \\ \hline
2       & 386,303   & 38,630        & 1,101,790     & 511,850       & 23.2 hours    \\ \hline
4       & 118,249   & 23,650        & 349,069       & 156,679       & 8.76 hours    \\ \hline
8       & 33,557    & 13,423        & 105,775       & 44,465        & 2.49 hours    \\ \hline
16      & 9,053     & 7,242         & 32,156        & 11,995        & 56.8 minutes  \\ \hline
32      & 2,283     & 3,652         & 9,935         & 3,025         & 11.2 minutes  \\ \hline
64      & 613       & 1,961         & 3,650         & 814           & 4.44 minutes  \\ \hline
128     & 157       & 1,000         & 1,446         & 215           & 1.31 minutes  \\ \hline
256     & 43        & 542           & 724           & 89.0          & 40.9 seconds  \\ \hline
512+4     & 8         & 197           & 481           & 142           & 24.8 seconds  \\ \hline
1024+24    & 1         & 52.4          & 1,107         & 527           & 1.20 minutes  \\ \hline
\hline

\end{tabular}
\end{center}
\label{Tab:TheoryScaling}

\end{table}

To see if a PDE matrix can also give results in line with the idealized Estimate 4.1, we refer back to the biharmonic problem in Example 3. 
While the reduction in matrix-vector products is not quite as large as in Example 11, it is substantial.  For example, going from no polynomial preconditioning to a degree 200 polynomial reduces the matrix-vector products by a factor of about 110.  This is over half of the ideal reduction factor of 200. 

\section{Comparison to FGMRES}

Here we compare polynomial preconditioned GMRES to the related method FGMRES and to a PP-GMRES variant where we change the polynomial at each cycle.

\subsection{Versus FGMRES}

The methods FGMRES~\cite{Sa93} and GMRESR~\cite{vdVVu} are related to PP-GMRES.  Both of these methods allow GMRES to have preconditioning that varies at every iteration.  If we choose this preconditioner to be a cycle of GMRES(d), this corresponds to using a new polynomial to precondition each iteration of the outer GMRES.  This is in contrast to PP-GMRES, which has a fixed polynomial that is used to precondition all outer GMRES iterations.  
These polynomials also differ in how they are generated: 
The polynomial for PP-GMRES is designed to approximate a solution for a linear equations problem with a random right-hand side.  This constrains it to have small norm over the entire spectrum of $A$ so that it can effectively precondition for any right-hand side.  
Meanwhile the GMRES(d) used for a step of FGMRES attacks the specific linear equations problem at that moment and thus may be skewed for that problem instead of applying a polynomial that addresses the overall spectrum of $A$.

In the comparisons that follow, we use the same degree polynomials for PP-GMRES and FGMRES.  FGMRES uses one more matrix-vector product per iteration because its GMRES$(d)$ has $d$ matrix-vector products in order to build a $p$ polynomial of degree $d-1$.  Then FGMRES still has to apply the operator $A$. 

{\it Example 12.} The matrix is diagonal with entries $\frac{1^2}{n}, \frac{2^2}{n},  \ldots, \frac{n^2}{n}$, the same matrix from Example 11, except here $n = 10{,}000$. Comparison between polynomial preconditioned GMRES and FGMRES is in Table \ref{Tab:FGMRESandChPoly}. (The third method in the table will be discussed in the next subsection.)  For low degree polynomials, FGMRES is significantly better than PP-GMRES in terms of both matrix-vector products and computation time.  This probably is because FGMRES continually changes its polynomial to be optimal for the current situation.  However, FGMRES requires orthogonalization to implement its GMRES cycle at every iteration, and so for high degree polynomials, PP-GMRES is much cheaper.  PP(d)-GMRES(50) reduces solve time to 5.7 seconds with $d = 250$, while the minimum solve time of FGMRES(50) is 253 seconds with $d = 200$.  Given a matrix with a more expensive matrix-vector product, the orthogonalization expense of FGMRES would use a smaller proportion of total solve time, and FGMRES might be the best method.  However, in a highly parallel setting, the global communication costs of orthogonalization become increasingly prohibitive so PP-GMRES may still have the advantage.

\begin{table}

\caption{Comparison with different degree polynomials between PP-GMRES (PP-G), FGMRES (FG) and PP-GMRES with polynomial changing for each cycle (ChPoly).  GMRES(50) is used for all tests.  The matrix is diagonal with entries $\frac{i^2}{n}$ and with $n=10{,}000$.}

\begin{center}
\begin{tabular}{|c||c|c|c||c|c|c|}  \hline\hline
        & PP-G          & FG            & ChPoly        & PP-G        & FG        & ChPoly \\
degree  & $mvps$        & $mvps$        & $mvps$        & time        & time      & time  \\
d       & (thousands)   & (thou's)      & (thou's)      & (seconds)   & (sec's)   & (sec's) \\ \hline \hline
5       & 5765          & 2468          & 580           & 3006        & 1166      & 191   \\ \hline
10      & 2798          & 1064          & 335           & 907         & 472       & 64.0  \\ \hline
25      & 1213          & 693           & 293           & 168         & 454       & 32.7  \\ \hline
50      & 593           & 383           & 365           & 39.6        & 386       & 36.7  \\ \hline
100     & 322           & 209           & 550           & 18.4        & 367       & 61.5  \\ \hline
150     & 223           & 155           & 636           & 12.0        & 389       & 88.9  \\ \hline
200     & 202           & 76.4          & 747           & 11.3        & 253       & 125   \\ \hline
250     & 95.4          & 82.5          & 626           & 5.7         & 338       & 128   \\ \hline
\hline

\end{tabular}
\end{center}
\label{Tab:FGMRESandChPoly}

\end{table}

\subsection{Polynomial preconditioning with the polynomial changed for each cycle}
\label{Sec:ChangePoly}

We now introduce a new variation of PP-GMRES where the polynomial is changed every time GMRES restarts.  We recompute the polynomial at the start of each GMRES cycle using the current residual vector as the starting vector for the polynomial computation.  
This is contrary to our advice in Subsection \ref{Sec:StartVec}; however the strategy works well in the following example.
We view this approach mainly as an interesting phenomenon, rather than as a competing method.  

{\it Example 12 (cont.)}  
Table \ref{Tab:FGMRESandChPoly} has columns labeled ChPoly where the polynomial is changed for each cycle of PP-GMRES.  This works surprisingly well for low degree polynomials, but becomes more expensive for high degree polynomials.  For low degree polynomials, the `changing polynomial' method needs far less iterations than FGMRES.  We suspect that this is because even though FGMRES changes polynomials more frequently, it is focused on solving a particular linear equations problem instead of the overall problem.  Even though the changing polynomial method struggles for high degree polynomials, it still stays ahead of FGMRES in solve time.  Ultimately its fastest time does not match that of regular PP-GMRES.  Changing the polynomial for PP-GMRES loses its advantage for high degrees, because if a residual vector is skewed at the beginning of a cycle, meaning some eigencomponents are much larger than others, then the polynomial for the preconditioning is skewed and less effective.  For low degree polynomials, a skewed polynomial is used for fewer matrix-vector products and is quickly replaced by another polynomial that can compensate for the previous one.  


\section{Double Polynomial Preconditioning}

Some earlier examples showed that polynomial preconditioning can greatly reduce dot products.  Here we look at further reducing dot products by using high degree composite polynomials.

In the results of Example 3 shown in Table \ref{Tab:BihamrN40Th}, dot products are reduced by almost four orders of magnitude going from no polynomial preconditioning to a polynomial of degree $200$.  However, there is a limit to this reduction; at some point, creating higher degree polynomials raises the total number of dot products.  To study this, we separate the two phases of PP-GMRES: polynomial creation and the linear solve.   When solving the linear equations, dot products keep going down as the cycles are reduced, but the dot products needed to generate the polynomial increase as the polynomial degree increases.  With degree $200$ in Example 3, there are about $20{,}000$ dot products for generating the polynomial and about $14{,}000$ for the linear solve.  Then with degree $400\ (+3)$, only about $5{,}000$ dot products are needed for the solve, but about $80{,}000$ dot products are needed for generating the polynomial.

For high-degree polynomials that reduce the total number of dot products,
we suggest double polynomial preconditioning~\cite{PPArnoldi}.  
First, a GMRES iteration for matrix $A$ of length $d_1$ finds the polynomial $\phi_1$.  Then GMRES with matrix $\phi_1(A)$ is run to length $d_2$ to determine the polynomial $\phi_2$.  The corresponding polynomials $p_1$ and $p_2$ are such that $\phi_1(\alpha) = \alpha p_1(\alpha)$ and $\phi_2(\alpha) = \alpha p_2(\alpha)$.    The composite polynomial, $\phi_2(\phi_1(A))$, is used for the polynomial preconditioned GMRES phase.  Plugging in to (\ref{eqn:ppsys}) and (\ref{eqn:ppsys2}), the linear equations problem becomes 
\begin{align*}
\phi_2(\phi_1(A))z &= b,  \\
x &= p_1(A)(p_2(\phi_1(A))z.  
\end{align*}
This approach allows for high degree polynomials with minimal storage and orthogonalization for the initial GMRES iteration which creates the composite polynomial.

{\it Example 13.}  We use the same matrix as in Example 3.  The results are in Table \ref{Tab:BiharmDouble}, which is somewhat a continuation of Table \ref{Tab:BihamrN40Th}.  
We use $d_1 = d_2 = 10, 15, 20, 30, 40, 50,$ and $60$ (for simplicity, the degrees of $\phi_1$ and $\phi_2$ are the same, but this is not necessary).  These values of $d_1$ and $d_2$ give very high degree composite polynomials, up to degree 3600.  The best results are for degree 1600 with the fastest time reduced from 47.9 seconds for a single degree 403 polynomial in Table \ref{Tab:BihamrN40Th} to 35.8 seconds here.  However, the more significant improvement is that dot products go down by a factor of more than 10 from the best single polynomial (degree 200) to the best composite polynomial (degree 1600).  Generating two polynomials of degree 40 takes many fewer dot products than generating one polynomial of degree 200.

\begin{table}

\caption{  Biharmonic matrix with $n=40,000$. }

\begin{center}
\begin{tabular}{|c|c|c|c|c|c|}  \hline\hline
degree & cycles    & $mvps$       & $vops$        & dot prod's  & time          \\
deg = $d_1$ x $d_2$  &           & (thousands)   & (thou's)   & (thou's)   &               \\  \hline \hline
100 = 10 x 10  & 117       & 583           & 903           & 154           & 5.29 minutes  \\ \hline
225 = 15 x 15  & 31        & 348           & 433           & 41.0          & 3.00 minutes  \\ \hline
400 = 20 x 20  & 8         & 153           & 174           & 10.2          & 1.26 minutes  \\ \hline
900 = 30 x 30  & 3         & 99.1          & 107           & 3.66          & 49.5 seconds  \\ \hline
1600 = 40 x 40  & 1         & 75.3          & 81.0          & 2.76          & 35.8 seconds  \\ \hline
2500 = 50 x 50  & 1         & 77.6          & 83.9          & 3.07          & 36.6 seconds  \\ \hline
3600 = 60 x 60  & 1         & 79.3          & 87.4          & 3.95          & 38.3 seconds  \\ \hline
\hline

\end{tabular}
\end{center}
\label{Tab:BiharmDouble}

\end{table}

The reduction in dot products is more remarkable when compared to unpreconditioned 
GMRES; the improvement is five orders of magnitude for this example.  As linear equations become larger and computer architectures necessitate low-communication algorithms, double polynomial preconditioning is one possible tool to create high-degree polynomials in a cost-effective manner.

\section{Conclusion}

Previous polynomial preconditioners for GMRES generally have had complicated implementations or instability at high degrees.  
In this paper we address both problems by presenting a new implementation of the GMRES polynomial.  It is cheaper and more stable than previous implementations of this polynomial.  Furthermore, it is simple: to run the new polynomial preconditioned GMRES algorithm, the user only needs to specify the degree of the polynomial and when to restart GMRES ($d$ and $m$). For some problems, the polynomial can greatly reduce computational costs compared to regular restarted GMRES, and this is more likely when the problem is difficult.  We show that polynomial preconditioning can especially reduce dot products, which may help to avoid expensive global communication in a parallel setting.  Polynomial preconditioning can effectively accelerate standard preconditioners such as ILU.  It should be considered for problems that are slow to converge even with such standard preconditioning.

The polynomials are adjusted with multiple added roots for stability control.  Examples show that for many problems, these added roots are not needed except with high degree polynomials; for other problems with outstanding eigenvalues, added roots are essential.  We also give a test to check whether the stability control is sufficient. 

Polynomial preconditioning works well for Example 1 even though the matrix is indefinite.  However, several problems can appear for the indefinite case.  One such difficulty is that even with a real spectrum, the polynomial may not have a minimum at the origin and, thus, the spectrum is still indefinite after polynomial preconditioning.  Also, adding roots for stability on one side of the spectrum may increase the volatility of the polynomial on the other side (see Example 8).   Future work will address indefinite matrices.  Possible solutions include damping the polynomial~\cite{LiXiVeYaSa,PPArnoldi} and shifting the operator used for generating the roots of the polynomial.

We also plan to apply this polynomial preconditioning to non-restarted methods, such as the conjugate gradient method for symmetric problems and BiCGStab and IDR for nonsymmetric problems.  These methods do not suffer slowing convergence due to restarting, but there is still great potential to reduce orthogonalization expense and improve stability for difficult indefinite and non-normal problems. 
We could also apply this polynomial preconditioner to the eigenvalue deflated method GMRES-DR~\cite{GMRES-E,GMRES-DR,PPG}.  In some sense, both polynomial preconditioning and eigenvalue deflation accomplish the same thing, so it would be interesting to analyze the differences between them and study cases where both are needed.

\section*{Acknowledgments} We appreciate the many helpful comments of the referees.  This study does not have any conflicts to disclose.

\bibliographystyle{siam}
\bibliography{morgan,embree,loe}

\end{document}